\documentclass[12pt]{article}

\usepackage{amsmath}
\usepackage{amssymb}
\usepackage{amsthm}

\newtheorem{thm}{Theorem}[section]
\newtheorem{prop}[thm]{Proposition}
\newtheorem{defn}[thm]{Definition}
\newtheorem{exam}[thm]{Example}
\newtheorem{Coro}[thm]{Corollary}
\newtheorem{Lemm}[thm]{Lemma}
\newtheorem{Rema}[thm]{Remark}
\newcommand{\mult}{\mathrm{mult}}
\newcommand{\ord}{\mathrm{ord}}
\newcommand{\Spec}{\mathrm{Spec}}
\newcommand{\emb}{\mathrm{emb}}

\begin{document}
\setcounter{section}{-1}

\begin{center}
{\Large Freeness of adjoint linear systems on threefolds
with terminal Gorenstein singularities or some quotient singularities}\\
\end{center}
\begin{center}
 Nobuyuki Kakimi
\end{center}
\begin{flushleft}
 Department of Mathematical Sciences, University of Tokyo, Komaba, Meguro, 
 Tokyo 153, Japan (e-mail:kakimi@ms318sun.ms.u-tokyo.ac.jp)  
\end{flushleft}
\begin{flushleft}
Abstract. 
We generalize the result of Kawamata concerning the strong version of
Fujita's freeness conjecture for smooth 3-folds to some singular cases, 
namely, Gorenstein terminal singularities, 
Gorenstein ${\mathbb Q}$-factorial terminal singularities and 
quotient singularities of type $1/r(1,1,1)$ and of type $1/r(1,1,-1)$.
We generalize furthermore the result of that 
to projective threefolds with only canonical singularities for 
canonical and not terminal singularities. 
It turns out that the estimates in the first three cases are better 
than the one for the smooth case, while it is not in the fourth case.
We also give explicit examples which show the estimate in the fourth case is
necessarily worse than the one for the smooth case. 
\end{flushleft}
\section{Introduction}
\hspace{1.5 em}
We recall related results which are previously known.
T.Fujita conjectured that, 
for a smooth projective variety $X$ and  an ample divisor $L$ on
$X$, the linear system $|K_X+mL|$ is free if $m \geq \dim X+1$.
Reider [Rdr] proved that this is the case if $\dim X=2$.
Ein and Lazarsfeld  [EL] proved that this is the case if $\dim X=3$.
Kawamata [K3] proved that this is the case if $\dim X=4$.
For a projective variety $X$ of dimension $2$ with some singularities,    
Ein and Lazarsfeld [EL] and Matsushita [M] extended the result of Reider [Rdr] to singular cases.
Kawachi [KM] obtained the effective estimates for a normal surface $X$.
For a projective variety $X$ of dimension $3$ with some singularities,  
Oguiso and Peternell [OP] proved that $|K_X +5L|$ is free
if $X$ is a projective threefold with only ${\mathbb Q}$-factorial terminal
Gorenstein singularities.
Ein, Lazarsfeld and Ma\c{s}ek [ELM] and Matsushita [M]
extended some of the results of Ein and Lazarsfeld [EL] 
to projective threefolds with terminal singularities.

A strong version of Fujita's freeness conjecture: \\
Let $X$ be a normal projective variety of dimension $n$, 
$x_0 \in X$ a smooth point, and $L$ an ample Cartier divisor.
Assume that there exist positive numbers $\sigma_p$ for $p=1,2, \cdots, n$
which satisfy the following conditions:
 (1) $\sqrt[p]{(L)^p \cdot W} \geq \sigma_p$ for 
any subvariety $W$ of dimension $p$ which contains $x_0$,
(2)$\sigma_p \geq n$ for all $p$, and $\sigma_n> n$.
Then $|K_X+L|$ is free at $x_0$.

Fujita [F] proved that, if $\sigma_1 \geq 3$, $\sigma_2 \geq \sqrt{7}$, 
and $\sigma_3 > \sqrt[3]{51}$, then $| K_X + L |$ is free at $x_0$. 
Our results are generalizations of the following result of Kawamata [K3]:
Let $X$ be a normal projective variety of dimension $3$,
$L$ an ample Cartier divisor, and $x_0 \in  X$ a smooth point.
Assume that there are positive numbers $\sigma_p$
for $p = 1,2,3$ which satisfy the following conditions:
(1) $\sqrt[p]{(L)^p \cdot W} \geq \sigma_p$ for 
any subvariety $W$ of dimension $p$ which contains $x_0$,
(2) $\sigma_1 \geq 3$, $\sigma_2 \geq 3$, and $\sigma_3 > 3$.
Then $| K_X + L |$ is free at $x_0$.

We shall prove the following results in this paper:
Let $X$ be a normal projective variety of dimension $3$,
$x_0 \in X$ a point, and $L$ an ample ${\mathbb Q}$-Cartier divisor
such that $K_X+L$ is a Cartier divisor at $x_0$.
Assume that $\sqrt[p]{(L)^p \cdot W} \geq \sigma_p$ for 
any subvariety $W$ of dimension $p$ which contains $x_0$.\\
$(1)$(Theorem 3.1)    
Let $x_0 \in X$ be a Gorenstein terminal singular point.
Assume that $\sigma_1 > (\sqrt[3]{2}+\sqrt{3})/\sqrt[3]{2}$,
$\sigma_2>\sqrt[3]{2}+\sqrt{3}$, and $\sigma_3>\sqrt[3]{2}+\sqrt{3}$. 
Then $| K_X + L | $ is free at $x_0$.
$(\mathrm{Note\ that\ }3>\sqrt[3]{2}+\sqrt{3}>2.99)$\\
$(2)$(Theorem 3.3)
Let $x_0 \in X$ be a Gorenstein terminal singular 
${\mathbb Q}$-factorial point.
Assume that $\sigma_1 \geq 2$, $\sigma_2 \geq 2 \sqrt{2}$, and 
$\sigma_3 >  2 \sqrt[3]{2}$.
Then $| K_X + L | $ is free at $x_0$.\\ 
$(3)$(Theorem 3.4) 
Let $x_0 \in X$ be a quotient singular point of type
$({1}/{r},{1}/{r},{1}/{r})$ for an integer $r$.
Assume that $\sigma_1 \geq {3}/{r}$, $\sigma_2 \geq {3}/{\sqrt{r}}$,
and $\sigma_3 > {3}/{\sqrt[3]{r}}$.
Then $| K_X + L |$ is free at $x_0$.\\
$(4)$(Theorem 3.6)
Let $x_0 \in X$ be a quotient singular point of type
$(1/r,1/r,-1/r)$ for an integer $r \geq 3$.
Assume that $\sigma_1 \geq 1+(1/r)$, $\sigma_2 \geq (1+(1/r))\sqrt{r+3}$, 
and $\sigma_3 > (1+(1/r))\sqrt[3]{r+2}$.
Then $| K_X + L |$ is free at $x_0$.\\
$(5)$(Theorem 3.8)
Let $X$ be a projective threefold with only canonical singularities and 
$x_0 \in X$ be a canonical and not terminal singular point.
Assume that $\sigma_1 \geq {3}$, $\sigma_2 \geq {3}$, and $\sigma_3 > {3}$.
Then $| K_X + L |$ is free at $x_0$.

The result (1), the result (5), and the result  
for the smooth case [K3] imply the following:
Corollary (Corollary 3.9). 
Let X be a projective variety of dimension $3$ 
and $H$ an ample Cartier divisor on $X$. Assume that $X$ has at most 
canonical Gorenstein singularities. Then $|K_X+mH |$ is free if $m \geq 4$.
Moreover, if $(H^3) \geq 2$, then $|K_X +3H|$ is also free. 
Note that Lee[L1, L2] also obtained a result similar to this corollary
independently. 

We also give a series of examples (Example 3.5) which shows that
our estimate in $(3)$ is optimal for each $r$.
Note that the estimates for $\sigma_p$ in cases $(1), (2)$ and 
$(3)$ are better than the one in the smooth case. On the other hand, 
our estimate in $(4)$ is worse than in the smooth case, especially
when the indices $r$ being large, indeed, 
$\sigma_3 \to +\infty$ if $r \to +\infty$.
In Example 3.7, we construct an infinite sequence $(X_r,L_r)
(r=1,2,3,\cdots)$ consisting of a $3$-fold $X_r$ having a singular point
of type $1/r(1,1,-1)$ and an ample ${\mathbb Q}$-Cartier divisor $L_r$ 
on $X_r$ which satisfies that $(L_r^3) \to +\infty$ if $r \to +\infty$, 
$K_{X_r}+L_r$ is Cartier at the point, 
but $|K_{X_r}+L_r|$ is not free at the point for each $r$.
This sequence shows that there exists no uniform estimate for $\sigma_3$ 
independent of the indices $r$ and also explains the reason 
why $\sigma_3 \to +\infty$ under $r \to +\infty$ in our estimate $(4)$.

Our proof is very similar to the one for the smooth case given in [K3].
However this involves more  careful and  detailed analysis of multiplicities
and discrepancies. For example, important differences for proof in between
the case $(1)$ and the smooth case are: the discrepancy coefficient for $K_X$ 
under the blow up $x_0 \in X$ is $1$ in $(1)$, 
while it is $2$ in the smooth case;
$\mult_{x_0}X=2$ in $(1)$, while it is $1$ in the smooth case
and the multiplicity $\mult_{x_0}S$ of 
the minimal center $S$ (Definition 1.2) at $x_0$ is $1$, $2$ or $3$ in $(1)$ 
while $\mult_{x_0}S= 1$ or $2$ in the smooth case if $S$ is a surface.
We also notice that in cases $(3)$ and $(4)$ the discrepancy coefficients 
for $K_X$ under the blow up $x_0 \in X$ are no more integers;
it is $3/r-1$ in $(3)$ and $1/r$ in $(4)$.
By this reason, we need to treat ${\mathbb Q}$-Cartier divisors
whose orders $\ord_{x_0}D$ at $x_0$ is of the form $d/r$ and 
especially for Theorem 3.4, the case $(3)$, we need Lemma 2.3, 
a generalization of [K3, Theorem 2.2] to ${\mathbb Q}$-Cartier divisors of
fractional orders at $x_0$. 

Acknowledgment:
The author would like to express his thanks to
Professor Yujiro Kawamata
for his advice and warm encouragement.
Also the author would like to express his thanks to
Referee for his useful and kind comments.    

 \section{Preliminaries}
\hspace{1.5 em}
Most of the results of this paper are the applications of 
the following vanishing theorem:
\begin{thm}[{[K1, V]}] 
Let $X$ be a smooth projective variety and 
$D$ a ${\mathbb Q}$-divisor.
Assume that $D$ is nef and big, and that the support
of the difference $ {}^{\lceil} {D} { }^{\rceil} - D $ is a normal crossing
divisor. Then $ H^p (X, K_X + {}^{\lceil} {D} { }^{\rceil}) = 0 $ 
for $p > 0$.
\end{thm}

We recall notation of [K3](cf [KMM]).
 \begin{defn} \normalfont
 Let $X$ be a normal variety and $D =\Sigma_i d_i D_i$ 
 an effective ${\mathbb Q}$-divisor such that $K_X + D$ is ${\mathbb Q}$-Cartier.
 If $\mu:Y \rightarrow X$ is an embedded resolution of the pair $(X,D)$,
 then we can write 
 \[ K_Y + F = {\mu}^* (K_X + D) \]
 with $F = {\mu}_{*}^{-1}D +\Sigma_j e_j E_j$ 
 for the exceptional divisors $E_j$.

 The pair $( X,D )$ is said to have only 
 \textit{log canonical singularities} $(LC)\\ 
 (\mathrm{resp. } \mathit{kawamata\ log\ terminal\ singularities} (KLT))$ 
 if $d_i \leq 1 (\mathrm{resp.} < 1)$
 for all $i$ and  $e_j \leq 1(\mathrm{resp.} <1)$ for all $j$.

 A subvariety $W$ of $X$ is said to be a \textit{center of
 log canonical singularities} for the pair $( X,D )$,
 if there is a birational morphism from a normal variety 
 $\mu:Y \rightarrow X$ and a prime divisor $E$ on $Y$ 
 with the coefficient $e \geq 1$
 such that ${\mu}(E) = W$. 
 The set of all the centers of log canonical singularities is denoted
 by $CLC( X,D )$.
 The union of all the subvarieties in $CLC(X,D)$ is denoted by
 $LLC(X,D)$ and called the \textit{locus of log canonical
 singularities} for $(X,D)$.
 For a point $ x_0 \in X$, we define
 $CLC ( X,x_0,D ) = \{ W \in CLC( X,D ) ; x_0 \in W \}$. 
 
 \end{defn} 

We shall use the following propositions
proved by Kawamata [K3].
Then we shall control the singularities of 
the minimal center of log canonical singularities
and replace the minimal center of log canonical
singularities by a smaller subvariety.
\begin{prop} [{[K3, 1.5,1.6]}]
Let $X$ be a normal variety and  $D$ an effective 
${\mathbb Q}$-Cartier divisor such that $K_X + D $ is ${\mathbb Q}$-Cartier.
Assume that  $X$ is $KLT$ and $(X,D)$ is $LC$.
If $W_1,W_2 \in CLC(X,D)$ and $W$ is 
an irreducible component of $W_1 \cap W_2$,
then $W \in CLC(X,D) $.
If $(X,D)$ is not  $KLT$ at a point $x_0 \in X $,
then there exists the unique minimal element $W_0$ of $CLC(X,x_0,D)$.
Moreover, $W_0$ is normal at $x_0$.
\end{prop}
\begin{prop}[{[K3, 1.9]}]  
Let $x_0 \in X$, $D$ and $W_0$ be as in Proposition $1.3$. 
Assume that $\dim W_0 = 2$. 
Then $W_0$ has at most a rational singularity at $x_0$.
Moreover, if $W_0$ is singular at $x_0$, and if $D'$ is 
an effective ${\mathbb Q}$-Cartier divisor on $X$ such that 
${\ord}_{x_0} D'|_{W_0} \geq 1$, then $\{ x_0 \} \in CLC(X,x_0,D+D')$.
\end{prop}
\begin{Rema} \normalfont
The result in Proposition $1.4$ is extended to higher dimensions ([K4]).
\end{Rema}

\begin{prop}[{[K3, 1.10]}] 
Let $x_0 \in X$, $D$ and $W_0$ be as in Proposition $1.3$.
Let $D_1$ and $D_2$ be effective ${\mathbb Q}$-Cartier divisors on $X$ 
whose supports do not contain $W_0$ and which  induce the same
${\mathbb Q}$-Cartier divisor on $W_0$. 
Assume that $(X,D+D_1)$ is $LC$ at $x_0$ and
there exists an element of $CLC(X,x_0,D+D_1)$ which is properly
contained in $W_0$. 
Then the similar statement holds for the pair $(X,D+D_2)$.
\end{prop}
   
\section{General method}
\hspace{1.5 em}
We can construct divisors which have high multiplicity at a given point
from the following lemma.
\begin{Lemm} 
Let $X$ be a normal and complete variety of dimension $n$,
$L$ a nef and big ${\mathbb Q}$-Cartier divisor, 
$x_0 \in X$ a point, and  $t$,$t_0$ rational numbers
such that $t>t_0>0$.
Then there exists an effective  ${\mathbb Q}$-Cartier divisor $D$ such that
$ D \sim_{\mathbb Q} t L $ and
\[ {\ord}_{x_0}D \geq (t_0+\epsilon )
 \sqrt[n]{ \frac{(L^n)}{ \mult_{x_0}X } } \]
which is a rational number for $0 \leq \epsilon 
\ll \sqrt[n]{ \mult_{x_0}X/ (L^n)}$.
\end{Lemm}

\begin{proof}
We take $r \in {\mathbb Q}$ such that $rL$ is Cartier and 
$t' \in {\mathbb Q}$ such that $t'=t/t_0 >1$.
By [K3, 2.1], there exists an effective ${\mathbb Q}$-Cartier 
divisor $rD'$ such that 
\[ {\ord}_{x_0}rD' \geq r \sqrt[n]{ \frac{(L^n)}{ \mult_{x_0}X } } .\] 
Therefore we may take $D=t_0 D'$.
\end{proof}
 
We generalize [K3, 2.3] in which a given point is assumed to be 
a Gorenstein $KLT$ point.
For our purpose, we need to consider the case 
where a given point is a $KLT$ point.
The following proposition is the key of the proofs of our main results.
\begin{prop}  
Let $X$ be a normal projective 
variety of dimension $n$,
$x_0 \in X$ a $KLT$ point,  
and $L$ an ample ${\mathbb Q}$-Cartier divisor such that 
$K_X + L$ is Cartier at $x_0$.
Assume that there exists an effective ${\mathbb Q}$-Cartier divisor $D$
which satisfies the following conditions: \\
$(1)$ $D \sim_{\mathbb Q} t L $ for a rational number $ t<1$,\\
$(2)$ $(X,D)$ is $LC$ at $x_0$,\\
$(3)$ $\{x_0\} \in CLC(X,D)$.\\
Then  $| K_X + L | $ is free at $x_0$.
\end{prop}

\begin{proof}
(cf [K3 2.3])
Let $r=\min\{s \in {\mathbb N}$; $sL$ is Cartier\}
and $D'$ be a general member of $| mrL |$ for 
$m \gg 0$ which passes through $x_0$.
Replacing $D$ by $(1- \epsilon_1)(D + \epsilon_2 D')$
for some $0 < \epsilon_i \ll 1/(rm)$, 
we may assume that $x_0$ is an isolated point of 
$LLC (X,D)$.
Let $\mu:Y \rightarrow X$ be an embedded resolution of the pair
$(X,D)$.
Then 
\[ K_Y + E + F_1+F_2  = {\mu}^{*} ( K_X + D ) \]   
where $E$ is a reduced divisor such that
$\mu( E ) = \{x_0\}$, $F_1$ is a divisor of the form 
$\sum_{j} f_{1j} F_{1j}$ 
such that $f_{1j} < 1$ and $x_0 \in \mu(F_{1j})$, 
and $F_2$ is a divisor of the form 
$\sum_{i} f_{2i} F_{2i}$ such that $x_0 \notin \mu(F_{2i})$.
Then
\[ K_Y + ( 1- t ){\mu}^{*} L \sim_{\mathbb Q}
{\mu}^{*}( K_X + L ) - E - F_1-F_2 .\]
Thus 
\[ H^1( Y,{\mu}^{*}( K_X + L )- E + 
{}^{\lceil}{- F_1}^{\rceil} -F_{2}')=0\]
where $F_{2}' = {\mu}^{*}( K_X + L )
-{}^{\lceil}{{\mu}^{*}( K_X + L )-F_2}^{\rceil}$ 
and we obtain a surjection
\[ H^0(Y,{\mu}^{*}(K_X + L) + {}^{\lceil} {-{F_1}}^{\rceil} -F_{2}')
\rightarrow
H^0(E,{\mu}^{*} (K_X + L)) \cong {\mathbb C}. \]
Since ${}^{\lceil} {-F_1}^{\rceil}-F_{2}'$ is effective and exceptional
over a neighborhood of $x_0$, 
\[ H^0( X, K_X + L ) \rightarrow H^0(E,{\mu}^{*}( K_X + L )) \]
is also surjective. 
\end{proof}

We generalize [K3, 2.2]
in which ${\mathbb Q}$-divisors have integral order at $x_0$.
For our purpose, we need to treat ${\mathbb Q}$-Cartier divisors 
of fractional orders at $x_0$.

\begin{Lemm}
Let $X$ be a normal projective variety of dimension $3$,
${x_0} \in X$ a quotient singular point of type
$({1}/{r},{1}/{r},{1}/{r})$ for an integer $r$,
$L$ an ample ${\mathbb Q}$-Cartier divisor such that $K_X + L$ is Cartier
at $x_0$,
$W$ a prime divisor with 
${\ord}_{x_0} W = {d}/{r} \geq {1}/{r}$ for an integer $d$,
and $e$,$k$ positive rational numbers such that $de \leq 1$ and
$({k}/{r})^3 < {(L)^3}/{r^2}$.
Assume that there exists an effective ${\mathbb Q}$-divisor $D$
such that  $D \sim_{\mathbb Q} L$ and ${\ord}_{x_0} D \geq {k}/{r}$, 
and moreover that $D \geq ekW$ for any such $D$.
Then there exists a real number $\lambda$ with $0 \leq \lambda < 1$
and $\lambda \leq max \{ 1-de ,(3de)^{-{1}/{2}} \}$
which satisfies the following condition:
if $k'$ is a positive rational number such that $k' > k$ and
\[ (\lambda \frac{k}{r} )^3 + ( \frac{1-de-\lambda }{ 1-\lambda } )^2
\{ (\frac{k'}{r} + \frac{\lambda de }{1-de-\lambda}\frac{k}{r} )^3
-(\frac{ \lambda k}{r} + 
\frac{\lambda de}{1-de-\lambda}\frac{k}{r} )^3 \} < \frac{(L)^3}{r^2},\]
then there exists an effective ${\mathbb Q}$-divisor $D$ such that 
$D \sim_{\mathbb Q} L$ and ${\ord}_{x_0} D \geq {k'}/{r}$.
(If $\lambda=1-de$,  
then the left hand side
of the above inequality should be taken as a limit.)
\end{Lemm}
\begin{proof}
(cf [K3, 2.2])
We have $\mult_{x_0} X=r^2$.
Let $\bar{k}$ = sup \{$q$; there exists an effective ${\mathbb
Q}$-divisor $D$ such that $D \sim_{\mathbb Q} L$ 
and ${\ord}_{x_0} D={q}/{r}$ \}.
Let us define a function  $\phi(q)$ for $q \in {\mathbb Q}$
with  $0 \leq q < \bar{k}$ to be the largest real number such that
$D \geq \phi(q) W$ whenever $D \geq 0$, $D \sim_{\mathbb Q} L$ and 
${\ord}_{x_0} D = {q}/{r}$.
Then $\phi$ is a convex function.
In fact, if ${\ord}_{x_0} D_i = {q_i}/{r}$ and $D_i 
=( \phi(q_i) + \epsilon_i) W +$ other components for 
$0 \leq \epsilon_i \ll 1$ and $i = 1$, $2$, then 
${\ord}_{x_0} (t D_1 + (1-t) D_2 )= t ({q_1}/{r})+(1-t) ({q_2}/{r})$ and
$tD_1 + (1-t)D_2 
= ( t(\phi(q_1) +\epsilon_1)+(1-t)(\phi(q_2)+\epsilon_2))W +$
other components,
hence $\phi(t q_1 +(1-t) q_2) \leq t \phi(q_1) +(1-t) \phi(q_2)$.
Since $\phi(k) \geq ek$, there exists a real number $\lambda$ such
that $0 \leq \lambda < 1$ 
and $\phi(q) \geq {e(q- \lambda k)}/({1-\lambda})$ for any $q$.

Let $m$ be a large and sufficiently divisible integer and 
$\nu : H^0 (X,mL)\rightarrow O_{X,{x_0}}(mL)
\cong O_{X,{x_0}}$ the evaluation homomorphism.
We consider subspaces $V_i ={\nu}^{-1}({m_{x_0}}^{i} )$
of $H^0 (X,mL)$ for integers $i$ 
such that ${\lambda k m}/{r} \leq i \leq {k'm}/{r}$.
First, we have 
\[ \mathrm{dim} V_{\lceil {\lambda km}/{r}\rceil}
\geq  \mathrm{dim} H^0 (X,mL) - r^2 \frac{ ( \lambda k m /r)^3}{3!} 
+ \mathrm{lower\ terms\ in\ } m. \]  
Let $D \in | mL |$ be a member corresponding to 
$h \in V_i$ for some $i$.
Since we have $D \geq \phi(ir/m) m W$,
the number of conditions in order for $h \in V_{i+1}$ is at most 
( the number of homogeneous polynomials of order 
$i-\phi(ir/m)m ({d}/{r})$ in $3$ variables )
$\times$ $\mult_{x_0} X$, i.e.,
\[ r^2 \frac{ ( i- \phi(ir/m) m d/r )^2}{2!}+
\mathrm{lower\ terms\ in\ } m. \]
Therefore, we have $k' < \bar{k}$ because
\[ r^2 \frac{ (\lambda k m/r )^3}{3!} +
\sum_{i=\lceil {\lambda km/r}\rceil}^
{k'm/r-1}
r^2 \frac{ ( i - \frac{ e(ir/m -\lambda
k)}{(1-\lambda)} m {d}/{r})^2}{2!} +
\mathrm{lower\ terms\ in\ } m \]
\[= r^2 \frac{ (\lambda k m/r )^3 }{3!} +
r^2 \frac{ ( \frac{1-de-\lambda}{1-\lambda} )^2 
\{ ( k'/r + \frac{ \lambda d e}{1-de-\lambda}{k}/{r})^3 -
( \lambda k/r +
\frac{ \lambda d e}{1-de-\lambda } {k}/{r})^3 \} m^3 }{3!}\]
\[+ \mathrm{lower\ terms\ in\ } m \]
\[< \frac{m^3 (L^3)}{3!} + \mathrm{lower\ terms\ in\ } m.\]
By [K3 2.2 last part], we have that $\lambda \leq {\max}
\{ 1-de,(3de)^{-{1}/{2}} \}$.
\end{proof}

The following theorem plays the important role in 
the proofs of our main results.

\begin{thm}[{[A, Corollary 6]}]
Let $S$ be a normal surface, $x_0$ be a point of $S$.
Suppose $S$ has a rational singularity at $x_0$.
Let $Z$ be the fundamental cycle.

Then 
\[ {\mult}_{x_0}S= -Z^2 \]
for all integers $k$,
\[ \dim m_{S,x_0}^k/ m_{S,x_0}^{k+1} 
= k {\mult}_{x_0} S + 1  \] 
where $m_{S,x_0}$ is the maximal ideal of $x_0$ in $S$ and  
\[ \emb {\dim}_{x_0} S = {\mult}_{x_0} S + 1 .\]
\end{thm}

The following lemma plays the supporting role in 
the proof of Theorem 3.6. 

\begin{Lemm}
Let $(X,x_0)$ be a quotient singularity 
${\mathbb C}^3/{\mathbb Z}_r(a/r,-a/r,1/r)$ 
with $0 < a < r$.
Let $S_i={\mathbb C}^2/{\mathbb Z}_r(i/r,1/r)$
for $i = a, -a$.
Then 
\[{\emb}{\dim}_{x_0} X = {\mult}_{x_0} X + 2  .\]
\[ {\mult}_{x_0} X 
={\mult}_{x_0} S_a +{\mult}_{x_0} S_{-a} .\]
\[ \emb {\dim}_{x_0} X
=\emb {\dim}_{x_0} S_a +\emb {\dim}_{x_0} S_{-a} .\]
\end{Lemm}
\begin{proof}
Let  
$(xy)^{w}{x^s}{z^u}, (xy)^{w}{y^t}{z^u}  \in m_{X,x_0}^{n}/
m_{X,x_0}^{n+1}$
where $m_{X,x_0}$ is the maximal ideal of $x_0$ in $X$.
Then
${x^s}{z^u}, {y^t}{z^u}  \in m_{X,x_0}^{n-w}/
m_{X,x_0}^{n-w+1}$.
We consider that 
${x^s}{z^u} \in m_{S_a,x_0}^{n-w}/
m_{S_a,x_0}^{n-w+1}$
where $m_{S_a,x_0}$ is the maximal ideal of $x_0$ in $S_a$
and that
${y^t}{z^u} \in m_{S_{-a},x_0}^{n-w}/
m_{S_{-a},x_0}^{n-w+1}$ 
where $m_{S_{-a},x_0}$ is the maximal ideal of $x_0$ in $S_{-a}$.
By Theorem 2.4., for $i= a, r-a$
\[ \dim m_{S_{i},x_0}^{n-w}/ 
m_{S_{i},x_0}^{n-w+1} 
= (n-w) {\mult}_{x_0} S_{i} + 1 . \]
$m_{{S_a},x_0}^{n-w}/ 
m_{S_{a},x_0}^{n-w+1} \cap 
m_{S_{-a},x_0}^{n-w}/ 
m_{S_{-a},x_0}^{n-w+1} = (z^r)^{n-w}$ for $0 \leq w < n$.
Then
\[ {\dim} m_{X,x_0}^{n}/m_{X,x_0}^{n+1} 
\]
\[= \sum_{w=0}^{n-1}
 \{ \dim m_{S_{a},x_0}^{n-w}/ m_{S_{a},x_0}^{n-w+1}
 +\dim m_{S_{-a},x_0}^{n-w}/ m_{S_{-a},x_0}^{n-w+1}-1 \}+1     
\]
\[= \sum_{w=0}^{n-1}
 \{ (n-w)({\mult}_{x_0} S_{a} + {\mult}_{x_0} S_{-a}) + 1 \}+1 \]
\[ = ({\mult}_{x_0} S_{a} + {\mult}_{x_0} S_{-a}) \frac{n^2}{2}
+({\mult}_{x_0} S_{a} + {\mult}_{x_0} S_{-a}+2) \frac{n}{2} +1\]  
Hence  \[ {\mult}_{x_0} X ={\mult}_{x_0} S_{a}+{\mult}_{x_0} S_{-a} .\]
\[ {\emb}{\dim}_{x_0} X ={\mult}_{x_0} S_{a} +{\mult}_{x_0} S_{-a} +2 .\]
\end{proof}

\section{Main Theorem}
\hspace{1.5 em}
First we consider Gorenstein terminal singular points.
\begin{thm} Let $X$ be a normal projective variety of dimension $3$,
$x_0 \in X$ a Gorenstein terminal singular point, and
$L$ an ample ${\mathbb Q}$-Cartier divisor such that 
$L$ is Cartier at $x_0$. 
Assume that there are positive numbers $\sigma_p$
for $p = 1,2,3$ which satisfy the following conditions: \\
$(1)$ $\sqrt[p]{(L)^p \cdot W} \geq \sigma_p$ for 
any subvariety $W$ of dimension $p$ which contains $x_0$,\\
$(2)$ $\sigma_1 > (\sqrt[3]{2}+\sqrt{3})/ \sqrt[3]{2}$, 
$\sigma_2 > \sqrt[3]{2}+\sqrt{3}$, 
and $\sigma_3 > \sqrt[3]{2}+\sqrt{3}$.\\
Then $| K_X + L |$ is free at $x_0$.
\end{thm}
\begin{proof}
Since $x_0$ is  a Gorenstein terminal singularity of dimension $3$, 
$x_0$ is an isolated hypersurface 
singularity of multiplicity 2 ([R, (1.1) Theorem]).
Let $U$ be a neighborhood at $x_0$.
Let $U \subset V$ be an embedding of $U$ as 
a hypersurface in a smooth fourfold $V$, and 
let $g:\widetilde{V} \longrightarrow V$
be the blowing-up of $V$ at $x_0$.
Let $\widetilde{U}$ be the proper transform of $U$ 
in $\widetilde{V}$, $f: \widetilde{U} \longrightarrow U$
the restriction of $g$ to $\widetilde{U}$, 
and $F \subset \widetilde{V}$  
the exceptional divisor of $g$.
Then
\[  K_{\widetilde{U}}= ( K_{\widetilde{V}} + \widetilde{U} ) |_{\widetilde{U}}
=(g^{*} K_V +3F + g^{*}U -2F )|_{\widetilde{U}}
=f^{*}K_U+ (F|_{\widetilde{U}}). \]
(cf.[ELM, Lemma 2.2 's proof, l.8--l.13])

Step 0. 
Let $t$ be a rational number such that 
$t>{2\sqrt[3]{2}}/{\sqrt[3]{(L^3)}}$.
Since $\sigma_3 > 2\sqrt[3]{2}$, we can take $t < 1$.
Let $t_0$ be a rational number such that 
$t_0=2\sqrt[3]{2}/ \sqrt[3]{{(L^3)}}-\epsilon$
for $0 \leq \epsilon \ll \sqrt[3]{ 2/(L^3)}$. 
By Lemma 2.1, there exists an effective ${\mathbb Q}$-Cartier
divisor $ D$
such that $ D \sim_{\mathbb Q} t L $ and 
${\ord}_{x_0} D \geq (t_0 +\epsilon) \sqrt[3]{(L^3)/2}$.
Hence ${\ord}_{x_0} D =2$.

Let $c$ be the 
\textit{log canonical threshold} of $ ( X,D )$ at $x_0$:
\[ c=\sup{ \{t \in {\mathbb Q}; \mbox{ $K_X +tD$ is $LC$ at $x_0$ } \} }. \]
Then $c \leq 1$.
Let $W$ be the minimal element of $CLC(X,x_0,cD)$.
If $W=\{x_0\}$, then $| K_X+L |$ is free at $x_0$ by Proposition 2.2, 
since $c t <1$.\\

Step 1. We consider the case in which $W =C$ is a curve.
By Proposition 1.3, $C$ is normal at $x_0$, i.e., smooth at $x_0$.
Since $ t < 1$, we have
$c t +(1-c)<1$. Since $\sigma_1 \geq 2$, there exists a rational number $t'$
with $c t +(1-c)<t'<1 $ and an effective ${\mathbb Q}$-Cartier divisor $D_C' $
on $C$ such that  $D_C' \sim_{\mathbb Q}(t'-c t )L|_C$ and 
${\ord}_{x_0}D_C' =2(1-c)$.
As in [K3, 3.1 Step1], there exists an effective ${\mathbb Q}$-Cartier 
divisor $D'$ on $X$ such that $D'\sim_{\mathbb Q}(t'-c t )L$ 
and $ D'|_C=D_C'$.
Let $D_1'$ be a general effective ${\mathbb Q}$-Cartier divisor 
on an affine neighborhood $U$ of $ x_0$ in $ X$ such that 
$D_1'|_{C \cap U} = D_C'|_{C \cap U}$ and ${\ord}_{x_0}D_1'=2(1-c)$.
Then we have  ${\ord}_{x_0}(cD+D_1')=2$, hence $\{x_0\} \in CLC(U,cD+D_1')$.
Let 
\[ c' =\sup{\{ t \in {\mathbb Q};
\mbox{ $K_X +(cD+tD_1')$ is $LC$ at $x_0$ } \} }. \]
Since $D_1'$ is chosen to be general, we have $c' > 0$.
We have an element $W'$ such that $W'\in CLC(X,x_0,cD+c'D_1')$ and
$W' \not\supseteq C$. By Proposition 1.3, $CLC(X,x_0,cD+c'D_1')$ has
an element which is properly contained in $C$.
By Proposition 1.6, we conclude that $(X,cD+c'D')$ is $LC$ at $ x_0$, and
$CLC(X,x_0,cD+c'D')$ has
an element which is properly contained in $ C$, i.e., $\{x_0\}$.\\

Step 2. We consider the case in which $W =S$ is a surface.
By Proposition 1.4, $S$ has at most a rational singularity at $x_0$.\\ 
Step2-1. We assume first that $S$ is smooth at $x_0$.
As in Step 1, we take a rational number $t'$, 
an effective ${\mathbb Q}$-Cartier divisor $D'$ on $X$
and a positive number $c'$
such that $c t +(1-c)<t'<1$, $D' \sim_{\mathbb Q}(t'-c t )L$, 
${\ord}_{x_0}D'|_S, =2(1-c)$, $(X,cD+c'D')$ is $LC$ at $x_0$, and 
that the minimal element $W'$ of $CLC(X,x_0,cD+c'D')$
is properly contained in $S$.
Thus we have the theorem when $W'=\{x_0\}$.

We consider the case in which $W'=C$ is a curve.
Since $t, t' < 1 $, we have $c t + c'(t'-c t)+(1-c)(1-c') <1$.
As in Step 1, we take a rational number $t''$, 
an effective ${\mathbb Q}$-Cartier divisor $D''$ on $X$
and a positive number $c''$ such that  
$c t + c'(t'-c t)+(1-c)(1-c') <t''<1$,  
$D''\sim_{\mathbb Q}(t''-c t-c'(t'-c t))L$, 
${\ord}_{x_0}D''|_C=2(1-c)(1-c')$,
$(X,cD+c'D'+c''D'')$ is $LC$ at $x_0$ and 
that the minimal element $W''$ of $CLC(X,x_0,cD+c'D'+c''D'')$ is
properly contained in $C$, i.e., $ \{ x_0 \} \in
CLC(X,x_0,cD+c'D'+c''D'')$.
  
Step 2-2. We assume that $S$ has a rational singularity at $x_0$.  
Since the embedding dimension of $X$ at $x_0$ is $4$ and
$x_0$ is also a singular point of $S$, 
the embedding dimension of $S$ at $x_0$ is $3$ or $4$.
Therefore $d: =\mult_{x_0}S= 2$ or $3$, because $x_0 \in S$ is 
a rational singular point([A, Corollary 6]).
Since $t>{2\sqrt[3]{2}}/{\sqrt[3]{(L^3)}}$ and 
$\sqrt[3]{(L^3)}> \sqrt[3]{2}+\sqrt{3}$
, we can take $t <{ 2 \sqrt[3]{2}/ (\sqrt[3]{2}+\sqrt{3})}$
, so $t/2 +\sqrt{3}/ \sigma_2< 1$.
Therefore, if $ c \geq {1}/{2}$, then 
\[ c t + \frac{ 2 \sqrt{d} (1-c)}{\sigma_2} < 1 .\]
In this case, we can take a rational number $t'$ and  
an effective ${\mathbb Q}$-Cartier divisor $D'$ on $X$
such that $ c t + 2\sqrt{d}(1-c)/\sigma_2 < t' < 1 $,
$D' \sim_{\mathbb Q}(t'-ct )L$ and 
${\ord}_{x_0}D'|_S = 2 ( 1 - c )$,
and proceed as in Step 2-1.

On the other hand, if  $ c \leq {1}/{2}$, then 
\[ c t + \frac{\sqrt{d}}{\sigma_2} < 1 .\]
We take $t'$ and $D'$ with $D' \sim_{\mathbb Q}(t'-ct ) L $ and 
${\ord}_{x_0}D'|_S = 1$, by Proposition 1.4,
$ \{ x_0 \} \in CLC( X,x_0,cD+D' )$.
As in Step 1, there exists $ c' $ such that $ 1 \geq c' > 0 $,
$ ( X,x_0,cD+c'D' )$ is $LC$ at $x_0$,
and that the minimal element of $W'$ of $CLC( X,x_0,cD+c'D' )$
is properly contained in $S$.
If $W'=\{x_0\}$ , we have the theorem.

We consider the case in which  $W'= C$ is a curve.
We have 
\[{\ord}_{x_0}( c D + c' D')|_S \geq 2 c + c'
= 2 - 2 ( 1 - c - \frac{c'}{2}). \]
Since $ \sqrt{d}/ \sigma_2  < t' - c t $, 
we can take $t'-ct<{\sqrt{d}/(\sqrt[3]{2}+\sqrt{3})}$. Then, 
\begin{eqnarray*}
& & c t + c'( t'-c t ) + \frac{2}{\sigma_1} 
( 1 - c - \frac{c'}{2}) \\
&<& c ({ \frac{2 \sqrt[3]{2}}{\sqrt[3]{2}+\sqrt{3}}}) 
+ \frac{2 \sqrt[3]{2}}{\sqrt[3]{2}+\sqrt{3}} ( 1 - c )
+ c' ({ \frac{ \sqrt{d}}{\sqrt[3]{2}+\sqrt{3}} 
-\frac{\sqrt[3]{2}}{\sqrt[3]{2}+\sqrt{3}}}) \leq 1
\end{eqnarray*}
The rest is the same as before.
\end{proof}

The following example shows
that the condition $\sigma_3 > 2 \sqrt[3]{2}$ in Theorem 3.1 is optimal.
\begin{exam}\normalfont
Let $X=\{ xy+z^2+t^2=0 \} \subset  {\mathbb P}^4$ and 
$x_0=(0:0:0:0:1)$.
Then $x_0$ is a Gorenstein terminal singular point and 
$K_X=\mathcal{O}(-3)$.
If $L=\mathcal{O}(3)$, then $|K_X+L|$ is free at $x_0$.
If $L=\mathcal{O}(2)$, then $|K_X+L|$ is not free at $x_0$ and
$L^3 X=16$.
\end{exam}

If we assume furthermore that $X$ is ${\mathbb Q}$-factorial at $x_0$,
we obtain better estimates for $\sigma_p$.
\begin{thm}
Let $X$ be a normal projective variety of dimension $3$, $L$ an
ample ${\mathbb Q}$-Cartier divisor, $x_0 \in X$ a Gorenstein terminal
${\mathbb Q}$-factorial singular point.
Assume that there are positive numbers ${\sigma}_p$ for $p = 1,2,3$ 
which satisfy the following conditions: \\
$(1)$ $\sqrt[p]{ (L)^p \cdot W } \geq {\sigma}_p $
for any subvariety $W$ of dimension $p$ which contains $x_0$,\\ 
$(2)$ ${\sigma}_1 \geq 2$, ${\sigma}_2 \geq 2\sqrt{2}$, 
and ${\sigma}_3 > 2\sqrt[3]{2}$.\\
Then $| K_X + L |$ is free at $x_0$.
\end{thm}
\begin{proof}
Since ${\sigma}_1 \geq 2$ 
and ${\sigma}_3 > 2\sqrt[3]{2}$,
Steps 0,1 are the same as Steps 0,1 of the proof of Theorem 3.1.\\
Step 2. We consider the case in which $W = S$ is a surface.
Since $X$ is ${\mathbb Q}$-factorial, terminal Gorenstein at $x_0$,
$X$ is factorial at $x_0$([K2, Lemma 5.1]).
In particular, $S$ is a Cartier divisor at $x_0$.
Then we have $2 > \ord_{x_0} cD \geq \ord_{x_0} S$.
Hence we have $\ord_{x_0}S =1$ and $\mult_{x_0} S = 2$.
There exists a rational number $t'$
with $ct + (1 - c) < t' < 1$.
By Lemma 2.1 and ${\sigma}_2 \geq 2 \sqrt{2}$, 
there exists an effective ${\mathbb Q}$-Cartier divisor $D_S$ 
on $S$ such that
$D_S \sim_{\mathbb Q} (t'-ct ) L|_S$ 
and ${\ord}_{x_0}D_S =2 (1 - c)$. \\
The rest is the same as Step 2-1 of the proof of Theorem 3.1.
\end{proof}

We consider non Gorenstein singular points in the following.

The following theorem is generalization to dimension of $3$
of the following result of Kawachi [KM]:
Let $S$ be a normal projective surface,
$x_0 \in S$ a quotient singular point of type $({1}/{r},{1}/{r})$
for an integer $r$ and $L$ a nef and big ${\mathbb Q}$-Cartier divisor 
such that $K_X + L$ is Cartier at $x_0$.
If $LC \geq {2}/{r}$ for any curve $C$ through $x_0$ 
and $\sqrt{L^2} > 2/\sqrt{r}$, then $| K_X + L |$ is free at $x_0$.\\
  
\begin{thm}
Let $X$ be a normal projective variety of dimension $3$,
$x_0 \in X$ a quotient singular point 
of type $({1}/{r},{1}/{r},{1}/{r})$
for an integer $r$, 
and $L$ an ample ${\mathbb Q}$-Cartier divisor such that 
$K_X + L$ is Cartier at $x_0$.
Assume that there are positive numbers ${\sigma}_p$ for $p = 1,2,3$
which satisfy the following conditions: \\
$(1)$ $\sqrt[p]{ (L)^p \cdot W} \geq {\sigma}_p$ 
for any subvariety $W$ of dimension $p$ which contains $x_0$,\\
$(2)$ ${\sigma}_1 \geq {3}/{r}$,
${\sigma}_2 \geq {3}/{\sqrt{r}}$, and 
${\sigma}_3 >  {3}/{\sqrt[3]{r}}$.\\
Then $| K_X + L |$ is free at $x_0$.
\end{thm}

\begin{proof}
Since $X={\mathbb C}^3/{\mathbb Z}_r (1/r,1/r,1/r)$ around $x_0$,
setting $\mathcal{O}_{{\mathbb C}^3,0}={\mathbb C } \{ x, y, z \}$, 
we have $\mathcal{O}_{X,x_0}
={ \mathcal{O}_{{\mathbb C}^3,0} }^{ {\mathbb Z}_r (1/r,1/r,1/r)}
={\mathbb C} \{ x^{a} y^{b} z^{c} | a +b + c =r, a \geq 0, b \geq 0, c \geq 0 \}$
and ${m}_{X,x_0}^n=(x^{a} y^{b} z^{c} | a + b +c =n r)$ 
where ${m}_{X,x_0}$ is the maximal ideal of $x_0$ in $X$. Therefore
\[ \dim_{\mathbb C} \mathcal{O}_{X,x_0}/{m}_{X,x_0}^n
= \sum_{k=0}^{n-1} {(kr+2)(kr+1)}/{2}=(r^2/3!)n^3 
+ \mathrm{ lower \ terms \ in \ } n.\] 
Hence ${\mult}_{x_0} X = r^2$.\\
Step 0. Let $t$ be a rational number such that  
$t >  ({3}/{\sqrt[3]{r}}) / \sqrt[3]{({L}^3)} $.
Since ${\sigma}_3 > {3}/{\sqrt[3]{r}}$, we can take $t < 1$.
Let $t_0$ be a rational number such that 
$ t_0 = ({3}/{\sqrt[3]{r}}) / \sqrt[3]{({L}^3)}-\epsilon$
for $0 \leq \epsilon \ll \sqrt[3]{r^2/(L^3)}$.
By Lemma 2.1, there exists an effective ${\mathbb Q}$-Cartier divisor
$D$ such that  $D \sim_{\mathbb Q} t L$ and 
${\ord}_{x_0} D \geq (t_0+\epsilon)\sqrt[3]{(L^3)/r^2}$.
Hence ${\ord}_{x_0}D =3/r$.\\
Let $f:{\bar{X}} \rightarrow X$ be the blowing-up of $X$
at $x_0$, $E \subset {\bar{X}}$ the exceptional divisor of $f$
, and ${\bar{D}}$ the proper transform of $D$ in ${\bar{X}}$.
Then 
\[ K_{\bar{X}} = {f}^{*}K_X + (3/r-1) E .\]
\[{f}^{*} D = \bar{D}+ (3/r) E . \] 
Let $c$ be the log canonical threshold of $ ( X,D )$ at $x_0$:
\[ c =\sup{ \{t \in {\mathbb Q}; \mbox{ $K_X +t D$ is $LC$ at $x_0$ } \} }.\]
Then $c \leq 1$.
Let $W$ be the minimal element of $CLC(X,x_0,cD)$.
If $W=\{x_0\}$, then $| K_X+L |$ is free at $x_0$ by Proposition 2.2, 
since $c t <1$.\\

Step 1. We consider the case in which $W =C$ is a curve.
By Proposition 1.3, $C$ is normal at $x_0$, i.e., smooth at $x_0$.
Since $\sigma_1 \geq 3/r$,
as in Step $1$ of the proof of Theorem $3.1$,  
we take a rational number $t'$, 
an effective ${\mathbb Q}$-Cartier divisor $D'$ on $X$
and a positive number $c'$ 
such that  
$c t+(1-c) < t'<1$, 
$D' \sim_{\mathbb Q}(t'-c t )L$,  
${\ord}_{x_0}D'|_C =(3/r)(1-c)$, 
$(X, c D +c' D')$ is $LC$ at $ x_0$, and
that the minimal element $W'$ of $CLC(X,x_0,cD+c'D')$ 
is properly contained in $ C, i.e., \{x_0\}$.\\
  
Step 2. We consider the case in which $W=S$ is a surface. 
Let $U$ be a neighborhood at $x_0$,
and $h: \widetilde{U} \longrightarrow U$ 
its closure of local universal cover, \\
i.e., $U={\mathbb C}^3/{\mathbb Z}_r(1/r,1/r,1/r)$
and $\widetilde{U} ={\mathbb C}^3$.
Let $\widetilde{S}$ be $h^{*} (S |_{U} )$,
and $x_1$ be $h^{-1} (x_0)$.
Since $x_0 \in S$ is a $KLT$ point [K3, 1.7] and 
$h$ is ramified only over $x_0$, $\widetilde{S} \subset {\mathbb C}^3$ is 
also irreducible and $x_1 \in \widetilde{S}$ is at most $KLT$.
Moreover $\widetilde{S} \subset {\mathbb C}^3$ is a hypersurface, 
$\widetilde{S}$ is also Gorenstein. Therefore, 
$x_1 \in \widetilde{S}$ is a rational Gorenstein point, that is, 
$x_1 \in \widetilde{S}$ is a smooth point or a rational double point.
If $x_1$ is a smooth point of $\widetilde{S}$, 
then $S |_{U} \cong {\mathbb C}^2/{\mathbb Z}_r(1/r,1/r)$.
Hence we have $\mult_{x_0} S = r$.
If $x_1$ is a rational double point of $\widetilde{S}$, 
that $\widetilde{S}$ is $A_{rn+1}$ type $(n \in {\mathbb Z}_{\geq 0})$,
i.e.,
$\widetilde{S} \cong \Spec{\mathbb C}[ x^{rn+2},y^{rn+2}, xy ]$, 
since $\widetilde{S}$ is invariant by 
the action of ${\mathbb Z}_r (1/r,1/r,1/r)$.
Since the action of ${\mathbb Z}_r$ on ${\mathbb C}^3$ is a scalar 
multiplication by $\xi$, we may assume that 
the equation of $x_1 \in \widetilde{S}$ is of the standard form 
under the same coordinates of ${\mathbb C}^3$.
Then by looking the action of ${\mathbb Z}_r$ on the equation and 
using the fact that the equation must be semi-invariant, 
we conclude the result. Since 
\[ S|_{U}=\Spec{\mathbb C}[
(x^{rn+2})^{r},(x^{rn+2})^{r-1}(x y),\ldots ,(x^{rn+2})(x y)^{r-1}, \]
\[(xy)^{r},(y^{rn+2})^{r}, \ldots, (y^{rn+2})(x y)^{r-1} ],\]
the embedding dimension of $S|_U$ is $2r +1$.
Hence we have ${\mult}_{x_0} S = 2r$ ([A, Corollary 6]).
Hence ${\mult}_{x_0} S = 2r$ and ${\ord}_{x_0} S = 2/r$
 or ${\mult}_{x_0} S =r$ and ${\ord}_{x_0} S = 1/r$.
Let $d:={\mult}_{x_0} S /r = 1$ or $2$.\\

Step 2-1.
We assume first that $d = 1$.
As in Step 1,
there exists a rational number $t'$
with $c t + (1 - c) < t' < 1$.
By Lemma 2.1 and ${\sigma}_2 \geq 3/ \sqrt{r}$,  
there exists an effective ${\mathbb Q}$-Cartier divisor $D_S$ 
on $S$ such that
$D_S \sim_{\mathbb Q} (t'-c t ) L|_S$ and 
${\ord}_{x_0}D_S =(3/r)(1 - c)$.\\
The rest of Step 2-1 is the same as Step 2-1 of the proof of 
Theorem 3.1.\\

Step 2-2.
We assume that $d = 2$.
As in Step 2-1,
we take a rational number $t'$ with $c t +\sqrt{2}(1-c) < t'$
and  an effective ${\mathbb Q}$-Cartier divisor $D'$ on $X$ with
$D' \sim_{\mathbb Q}(t'-c t)L$ and
${\ord}_{x_0}D'|_S =(3/r)(1-c)$.
Here we need the factor $\sqrt{2}$
because $S$ has multiplicity $2r$ at $x_0$.
Then we take $0 < c' \leq 1$ such that 
$(X, c D +c' D')$ is $LC$ and $CLC(X,x_0,cD+c'D')$ has 
an element which is properly contained in $S$.

We shall prove that we may assume $c t +\sqrt{2}(1-c) < 1$.
Then we can take $t' < 1$ as in Step 2-1,
and the rest of the proof is the same.
For this purpose, we apply Lemma 2.3.
In argument of Steps 0 through 2-1, the number $t$ was chosen
under the only condition that $t<1$.
So we can take $t=1-\epsilon_1$, where the $\epsilon_n$ for 
$n = 1,2,...$ will stand for very small positive rational numbers.
Then $k/r= 3/(r(1-\epsilon_1)) 
=(3/r) +\epsilon_2$ and $e=1/(3c)$.
This means the following:
for any effective $D \sim_{\mathbb Q} tL$, 
if ${\ord}_{x_0} D \geq 3/r$,
then $cD \geq S$.
We look for $k'> 6/(3-\sqrt{2})$
so that there exists an effective ${\mathbb Q}$-Cartier divisor
$D \sim_{\mathbb Q} tL$ with $t < (3-\sqrt{2})/2$
and ${\ord}_{x_0} D \geq 3/r$.
The equation for $k'$ becomes
\[ {\lambda}^{3} +( \frac{1-2e-\lambda}{1-\lambda})^2
\{(\frac{2}{3-\sqrt{2}}+\frac{2 \lambda e}{1-2e-\lambda})^{3}
-(\lambda + \frac{2 \lambda e}{1-2e-\lambda})^{3} \} < 1.  \]
We have $\lambda \leq 1/ \sqrt{6e}$, 
$1/3 \leq e \leq 1/2$, and 
in particular, $0 \leq \lambda \leq 1/ \sqrt{2}$.
By [K3, 3.1 Step2-2],
we obtain a desired $D$, and can choose a new $t$
such that $t < (3-\sqrt{2})/2$.
Then we repeat the preceding argument from Step 0.
If we arrive at Step 2-2 again, then we have 
$2/3 \leq c \leq 1$ and $ct +\sqrt{2}(1-c) < 1$.  
\end{proof}

The following example shows that the conditions in Theorem 3.4
is best possible.

\begin{exam}\normalfont
Let $X={\mathbb P}(1,1,1,r)$ and $x_0=(0:0:0:1)$.
Then $x_0$ is a quotient singular point of type $(1/r,1/r,1/r)$ and 
$K_X = \mathcal{O}(-3-r)$.
If $K_X + L$ is Cartier and $L$ is effective,
we have $L = \mathcal{O}(rk+3)$ $(k \in {\mathbb Z},rk+3 \geq 0)$.\\
If $L = \mathcal{O}(3)$,
then $| K_X + L |$ is not free at $x_0$.
Hence the condition ${\sigma}_3 > 3 /\sqrt[3]{r}$ 
is necessary in Theorem $3.4$.
\end{exam}

We consider non Gorenstein terminal singular points in the following. 
\begin{thm}
Let $X$ be a normal projective variety of dimension $3$,
$x_0 \in X$ a quotient singular point
of type $({1}/{r},{1}/{r},{-1}/{r})$
for an integer $r \geq 3$, 
and $L$ an ample ${\mathbb Q}$-Cartier divisor such that 
$K_X + L$ is Cartier at $x_0$.
Assume that there are positive numbers ${\sigma}_p$ for $p = 1,2,3$
which satisfy the following conditions: \\
$(1)$ $\sqrt[p]{ (L)^p \cdot W} \geq {\sigma}_p$ 
for any subvariety $W$ of dimension $p$ which contains $x_0$,\\
$(2)$ ${\sigma}_1 \geq 1+(1/r)$,
${\sigma}_2 \geq (1+(1/r)){\sqrt{r+3}}$, and 
${\sigma}_3 >  (1+(1/r)){\sqrt[3]{r+2}}$.\\
Then $| K_X + L |$ is free at $x_0$.
\end{thm}

\begin{proof}
Since 
\[{\mathbb C}^3/{\mathbb Z}_r(1/r,1/r,-1/r) =\Spec{\mathbb C}[
x^{r}, x^{r-1}y ,\ldots, x y^{r-1},y^{r}, x z, y z, z^{r}],\] 
by Lemma 2.5, ${\mult}_{x_0} X = {\emb}{\dim}_{x_0} X - 2 = r + 2$.
Hence ${\mult}_{x_0} X = r+2$.
Step 0. Let $t$ be a rational number such that  
$t >  (1+(1/r))\sqrt[3]{r+2} / \sqrt[3]{({L}^3)}$.
Since ${\sigma}_3 > (1+(1/r))\sqrt[3]{r+2}$, we can take $t < 1$.
Let $t_0$ be a rational number such that 
$t_0 = (1+(1/r))\sqrt[3]{r+2} / \sqrt[3]{({L}^3)}-\epsilon$
for $0 \leq \epsilon \ll  \sqrt[3]{(r+2)/(L^3)}$.
By Lemma 2.1, there exists an effective ${\mathbb Q}$-Cartier divisor
$D$ such that  $D \sim_{\mathbb Q} t L$ and 
${\ord}_{x_0} D  \geq (t_0+\epsilon)\sqrt[3]{(L^3)/(r+2)}$.
Hence ${\ord}_{x_0}D=1+(1/r)$.\\
Let $f:{\bar{X}} \rightarrow X$ be the weighted blowing-up of $X$
at $x_0$, $E \subset {\bar{X}}$ the exceptional divisor of $f$
, and ${\bar{D}}$ be the proper transform of $D$ in ${\bar{X}}$.
Then 
\[ K_{\bar{X}} = {f}^{*}K_X + (1/r) E .\]
\[{f}^{*} D = \bar{D}+ e E , e \geq  1+(1/r). \]
Let $c$ be the log canonical threshold of $ ( X,D )$ at $x_0$:
\[ c =\sup{ \{t \in {\mathbb Q}; \mbox{ $K_X +t D$ is $LC$ at $x_0$ } \} }.\]
Then $c \leq 1$.
Let $W$ be the minimal element of $CLC(X,x_0,cD)$.
If $W=\{x_0\}$, then $| K_X+L |$ is free at $x_0$ by Proposition 2.2, 
since $c t <1$.\\

Step 1. We consider the case in which $W =C$ is a curve.
By Proposition 1.3, $C$ is normal at $x_0$, i.e., smooth at $x_0$.
Since $\sigma_1 \geq 1+(1/r)$, 
as in Step $1$ of the proof of Theorem $3.1$, 
we take a rational number $t'$,
an effective ${\mathbb Q}$-Cartier divisor $D'$ on $X$ 
and a positive number $c'$
such that  
$c t+(1-c)<t'<1$,
$D' \sim_{\mathbb Q}(t'-ct )L$,  
${\ord}_{x_0}D'|_C =( 1+(1/r) )(1-c)$,
$(X,cD+c'D')$ is $LC$ at $ x_0$, and
that the minimal element $W'$ of $CLC(X,x_0,cD+c'D')$ is 
properly contained in $ C, i.e., \{x_0\}$.\\
  
Step 2. We consider the case in which $W=S$ is a surface. 
Since the embedding dimension of $X$ at ${x_0}=r+4$, 
the embedding dimension of $S$ at ${x_0} \leq r+4$.
Hence we have ${\mult}_{x_0} S \leq r+3$ ([A, Corollary 6]).
There exists a rational number $t'$
with $c t + (1 - c) < t' < 1$.
By Lemma 2.1 and ${\sigma}_2 \geq (1+(1/r))\sqrt{r+3}$, 
there exists an effective ${\mathbb Q}$-Cartier
$D_S$ on $S$ such that
$D_S \sim_{\mathbb Q}(t'-ct ) L|_S$ 
and ${\ord}_{x_0}D_S =(1+(1/r))(1 - c)$. \\
The rest of Step 2 is the same as Step 2-1 of the proof of Theorem 3.1.
\end{proof}

The following example shows that the estimates in some terminal 
singular cases is necessarily worse than
the one for the smooth case.
\begin{exam}\normalfont
Let $X={\mathbb P}(1,1,r-1,r)$ and $x_0=(0:0:0:1)$.
Then $x_0$ is a quotient singular point of type $(1/r,1/r,-1/r)$ and
$K_X = \mathcal{O}(-2r-1)$.
If $K_X+L$ is Cartier at $x_0$ and $L$ is effective,
we have $L =\mathcal{O}(rk+1)$ $(k \in {\mathbb Z},rk+1 \geq 0)$.
If $L =\mathcal{O}(r+1)$,
then $|K_X + L|$ is not free at $x_0$ and 
$L^3 = (r+1)^3/r(r-1)$.
Hence if $r \geq 23$, then the condition ${\sigma}_3 >3$.  
\end{exam}

We obtain the following theorem based on [K3 Theorem 3.1] 
and [R, Main Theorem].
\begin{thm} Let $X$ be a projective variety of dimension $3$
with only canonical singularities,
$x_0 \in X$ be a canonical and not terminal singularity point, 
$L$ an ample ${\mathbb Q}$-Cartier divisor on $X$ 
such that $K_X+L$ is Cartier at $x_0$. 
Assume that there are positive numbers $\sigma_p$
for $p = 1,2,3$ which satisfy the following conditions: \\
$(1)$ $\sqrt[p]{(L)^p \cdot W} \geq \sigma_p$ for 
any subvariety $W$ of dimension $p$ which contains $x_0$,\\
$(2)$ $\sigma_1 \geq 3$, $\sigma_2 \geq 3$, and $\sigma_3 > 3$.\\
Then $| K_X + L |$ is free at $x_0$.
\end{thm}

\begin{proof}
By a theorem of Reid [R, Main Theorem], 
a partial resolution $f: Y \rightarrow X$
such that $K_Y =f^{*}K_X$, and
$Y$ has only terminal singularity points.
There exists a smooth point $y_0 \in f^{-1}(x_0)$.\\
Step 0.
Let $t$ be a rational number such that  
$t >  3 / \sqrt[3]{({f^{*}L}^3)} $.
Since ${\sigma}_3 > {3}$, we can take $t < 1$.
Let $t_0$ be a rational number such that 
$ t_0 = {3}/ \sqrt[3]{({f^{*}L}^3)}-\epsilon$
for $0 \leq \epsilon \ll \sqrt[3]{1/(L^3)}$.
By Lemma 2.1, there exists an effective ${\mathbb Q}$-Cartier divisor
$D$ such that  $D \sim_{\mathbb Q} t L$ and 
${\ord }_{y_0} f^{*}D \geq (t_0+\epsilon)\sqrt[3]{(L^3)}$.
Hence ${\ord}_{y_0}f^{*}D =3$.\\
Let $g:{\bar{Y}} \rightarrow Y$ be the blowing-up of $Y$
at $y_0$, $E \subset {\bar{Y}}$ the exceptional divisor of $g$
, and ${\bar{f^{*}D}}$ the proper transform of $f^{*}D$ in ${\bar{Y}}$.
Then 
\[ K_{\bar{Y}} = {g}^{*}K_Y + 2 E .\]
\[g^{*}f^{*} D = \bar{f^{*}D}+ 3 E . \] 
\[ K_{\bar{Y}} +\bar{f^{*}D}+ E = {g}^{*}f^{*}(K_X +D) .\]
Let $c$ be the log canonical threshold of $ ( X,D )$ at $x_0$:
\[ c=\sup{ \{t \in {\mathbb Q}; \mbox{ $K_X +t D$ is $LC$ at $x_0$ } \} }.\]
Then $c \leq 1$.
Let $W$ be the minimal element of $CLC(X,x_0,cD)$.
If $W=\{x_0\}$, then $| K_X+L |$ is free at $x_0$ by Proposition 2.2, 
since $c t <1$.\\

Step 1. We consider the case in which $W =C$ is a curve.
By Proposition 1.3, $C$ is normal at $x_0$, i.e., smooth at $x_0$.
Since $\sigma_1 \geq 3$,
as in Step $1$ of the proof of Theorem $3.1$,  
we take a rational number $t'$, 
an effective ${\Bbb Q}$-Cartier divisor $D'$ on $X$
and a positive number $c'$ 
such that  
$c t+(1-c) < t'<1$, 
$D' \sim_{\mathbb Q}(t'-c t )L$,  
${\ord}_{x_0}D'|_C =3(1-c)$, 
$(X, c D + c' D')$ is $LC$ at $ x_0$, and
that the minimal element $W'$ of $CLC(X,x_0,cD+c'D')$ 
is properly contained in $ C, i.e., \{x_0\}$.\\

Step 2. We consider the case in which $W=S$ is a surface.
We have $S' =f^{-1} S \in CLC(Y,y_0,f^{*}c D)$.
By Proposition 1.4,
$S'$ has at most a rational singularity at $y_0$. 
Step 2-1.
We assume first that $S'$ is smooth at $y_0$.
As in Step 1,
there exists a rational number $t'$
with $c t + (1 - c) < t' < 1$.
By Proposition 2.1 and ${\sigma}_2 \geq 3$,  
there exists an effective ${\mathbb Q}$-Cartier divisor $f^{*}(D_S)$ 
on $S'$ such that
$D_S \sim_{\mathbb Q} (t'-c t ) L|_S$ and 
${\ord}_{y_0}f^{*}(D_S) =3(1 - c)$.\\
The rest of Step 2-1 is the same as Step 2-1 of the proof of 
Theorem 3.1.\\
Step 2-2.
We assume that $d = 2$.
As in Step 2-1,
we take a rational number $t'$ with $c t +\sqrt{2}(1-c) < t'$
and  an effective ${\mathbb Q}$-Cartier divisor $D'$ on $X$ with
$D' \sim_{\mathbb Q}(t'-c t)L$ and
${\ord}_{y_0}f^{*}(D'|_S) =3(1-c)$.
Here we need the factor $\sqrt{2}$
because $S'$ has multiplicity $2$ at $y_0$.
Then we take $0 < c' \leq 1$ such that 
$(X, c D +c' D')$ is $LC$ and $CLC(X,x_0,cD+c'D')$ has 
an element which is properly contained in $S$.

We shall prove that we may assume $c t +\sqrt{2}(1-c) < 1$.
Then we can take $t' < 1$ as in Step 2-1,
and the rest of the proof is the same.
For this purpose, we apply Lemma 2.3.
In argument of Steps 0 through 2-1, the number $t$ was chosen
under the only condition that $t<1$.
So we can take $t=1-\epsilon_1$, where the $\epsilon_n$ for 
$n = 1,2,...$ will stand for very small positive rational numbers.
Then $k= 3/(1-\epsilon_1) 
=3 +\epsilon_2$ and $e=1/(3c)$.
This means the following:
for any effective $D \sim_{\mathbb Q} t L$, 
if ${\ord}_{y_0} f^{*}D \geq 3$,
then $c f^{*}D \geq S'$.
We look for $k'> 6/(3-\sqrt{2})$
so that there exists an effective ${\mathbb Q}$-Cartier divisor
$D \sim_{\mathbb Q} t L$ with $t < (3-\sqrt{2})/2$
and ${\ord}_{y_0} f^{*}D \geq 3$.
The equation for $k'$ becomes
\[ {\lambda}^{3} +( \frac{1-2e-\lambda}{1-\lambda})^2
\{(\frac{2}{3-\sqrt{2}}+\frac{2 \lambda e}{1-2e-\lambda})^{3}
-(\lambda + \frac{2 \lambda e}{1-2e-\lambda})^{3} \} < 1.  \]
We have $\lambda \leq 1/ \sqrt{6e}$, 
$1/3 \leq e \leq 1/2$, and 
in particular, $0 \leq \lambda \leq 1/ \sqrt{2}$.
By [K2, 3.1 Step2-2],
we obtain a desired $f^{*}D$, and can choose a new $t$
such that $t < (3-\sqrt{2})/2$.
Then we repeat the preceding argument from Step 0.
If we arrive at Step 2-2 again, then we have 
$2/3 \leq c \leq 1$ and $c t +\sqrt{2}(1-c) < 1$.  
\end{proof}
  
Theorem 3.1, Theorem 3.8, and 
the result for the smooth case [K3 Theorem 3.1]
imply the following corollary in which the estimate is better 
than the one in [OP, Theorem 2]
 
\begin{Coro}
Let $X$ be a projective variety of dimension $3$ 
with only \\Gorenstein canonical singularities, and
$H$ an ample Cartier divisor. \\
Then $| K_X + m H |$ is free if $m \geq 4$.
Moreover, if $(H^3) \geq 2$, then 
$| K_X + 3H |$ is also free.
\end{Coro}

We obtain the following corollary from Corollary 3.9 and [OP Theorem 1,3] 
\begin{Coro}[{ cf [OP Theorem I,II]}]
 Let $(X,L)$ be a polarized canonical Calabi-Yau threefolds.
$(1)$ $| m L |$ gives a birational map when  $m \geq 5$.
$(2)$ $| m L |$ is free if $m \geq 4$.
$(3)$ $| m L |$ is very ample when $m \geq 10$.
\end{Coro}

\end{document}